UDC 511.3

# The proof of the correctness of the Birch and Swinnerton–Diyer conjecture.

S.V.Matnyak

Khmelnitsky, Ukraine

**Abstract**. The proof of the conjecture of the Birch and Swinnerton-Dyer is presented in the *paper*. The Riemann's hypothesis on the distribution of non-trivial zeroes of the zeta-function of Riemann, previously proven, is word to prove this hypothesis. The theorem proved about the behavior of the $L$-function curve $E$ for $s \to 1$. It is shown that the $L$-function of the curve $E$ tends to zero for any prime unpaired integers. It is shown that the function can be expanded in a power series of the holomorphic field. The theorem is proved on conformity of the basis of the Galois group and the number of zero coefficients of the power series. The result proved the conjecture of Birch and Swinnerton-Dyer.

**Introduction**. In the beginning of 1960s, Birch and Swinnerton-Dyer set that dimension $r$ of the group of an elliptic curve $E$ over $Q$ is equal to the order of zero of Hasse – Weil zeta function $L(E,s)$ in the point $s=1$.

A more detailed conjecture states that there is no zero border set $B_E = \lim\limits_{s \to 1} \dfrac{L(E,s)}{(s-1)^r}$, where value $B_E$ depends on refined arithmetic invariants of curves.

The most important partial result for 2011 is be statement in 1977 by John Coates and Andrew Wiles state that is true for a large class of elliptic curves, which cays that if the curve $F$ has infinite amount of rational points, than $L(E,s)=0$.

The conjecture is the only relatively simple overall method of calculating the rank of elliptic curves.

1. **The problem (conjecture).** Let us consider that $E$ is a certain elliptic curve, stated off $Q$. Then the group rank $E$, $r_E$ is equal to zero order of $L$ - function $L(E,s)$, in point $s=1$.

   **Solution.** Let us assume that $E$ is the elliptic curve, stated off $Q$ with equation
   $$x_0 x_2^2 = x_1^3 - A \cdot x_0^2 \cdot x_1 - B \cdot x_0^3, \qquad A, B \in Q. \tag{1}$$

We obtain affine equation, when we set $x = \dfrac{x_1}{x_0}$ and $y = \dfrac{x_2}{x_0}$:

$$y^2 = x^3 - A \cdot x - B. \tag{2}$$

The transformation $(x, y) \to (c^2 x, c^2 y)$ turns this equation into

$$y^2 = x^3 - c^4 A x - c^6 B. \tag{3}$$

So, from the beginning we can assume that $A, B \in Z$. The Number $\Delta = 16(4 \cdot A^2 - 27B^2)$ is denominated as discriminant of curve $E$. As we can see $\Delta \neq 0$. Let the $p \in Z$ be some prime number, then we consider the equation

$$y^2 = x^3 - Ax - B(p),$$

That is the equivalent to the equation

$$y^2 = x^3 - \bar{A}x - \bar{B}, \quad \bar{A}, \bar{B} \in Z/pZ = F_p. \tag{4}$$

This equation determines the elliptical curve $E_p$ off $F_p$, only when $(p, \Delta) = 1$. Later will be considered only such primes, when the approved is not determined. The curve $E_p$ is called the reduction of curve $E$ on the module $p$.

Let the $N_{p^m}$ denotes the amount of points in $E_p(F_{p^m})$. Than we can analyze zeta function

$$Z(E_p, u) = \exp\left(\sum_{m=1}^{\infty} N_{p^m} \frac{u^m}{m}\right). \tag{5}$$

Using the Riemann-Roch theorem, we can write

$$Z(E_p, u) = \frac{1 - a_p u + pu^m}{(1-u) \cdot (1-pu)}, \quad a_p \in Z. \tag{6}$$

We write down for $a_p^2 \leq 4p$

$$1 - a_p u + pu^2 = (1 - \pi \cdot u)(1 - \bar{\pi} \cdot u), \tag{7}$$

Where $\bar{\pi}$ is complexly conjugated with $\pi$. It is clear that $\pi \cdot \bar{\pi} = p$, $a_p = \pi + \bar{\pi}$.

Besides, $|\pi| = |\bar{\pi}| = \sqrt{p}$. This is Riemann's conjecture for elliptical curve off $F_p$. Logarithmically differentiating (5) and (6), and taking into account (7) and comparing the coefficients, we will have

$$N_{p^m} = p^m + 1 - \pi^m - \bar{\pi}^m. \tag{8}$$

Especially, $N_p = p + 1 - a_p$. In such way calculating $N_p$, we determine $a_p$. Considering $\pi$ and $\bar{\pi}$ are roots of the equation $T^2 - a_p \cdot T + p = 0$, than the equation (8) determines $N_{p^m}$ for all $m \geq 1$.

We change variable $u$ for $p^{-s}$ and will get

$$\varsigma(E_p, s) = \frac{1 - a_p \cdot p^{-s} + p^{1-s}}{(1 - p^{-s}) \cdot (1 - p^{1-s})}, \tag{9}$$

We have determined $\varsigma(E_p, s)$ for prime numbers $(p, \Delta) = 1$. When $p/\Delta$, we consider

$$\varsigma(E_p, s) = \frac{1}{(1-p^{-s}) \cdot (1-p^{1-s})}.$$

Now, having introduced in the local zeta function for all primes p, we determine the global zeta function simply as the product of local zeta functions

$$\varsigma(E, s) = \prod_3 \varsigma(E_p, s). \tag{10}$$

From the definition we can see that,

$$L(E, s) = \frac{\varsigma(s) \cdot \varsigma(1-s)}{\varsigma(E, s)}. \tag{11}$$

We write down the function (11) as

$$L(E, s) = \frac{\varsigma(s) \cdot \varsigma(1-s)}{\varsigma(E, s)} = \frac{-0.5 \cdot k_1 \cdot \prod_p (1-p^{-s}) \cdot (1-p^{-1-s})}{\prod_p (1 - a_p p^{-s} + p^{1-2s})}, \tag{12}$$

where $k_1$ is the value of function $\frac{1}{\varsigma(1)}$. The function $\varsigma(1)$ coincides with the theorem 2 [3, p.2]

The theorem 1. The function $L(E, s)$ at $s = 1$ will be equal to zero $(L(E, s) = 0)$ with all values of $p$.

**Proof.** We determine

$$k_1 = \frac{1}{\varsigma(1)} = \sum_{n=1}^{\infty} M(n) \cdot \left(\frac{1}{n+1} - \frac{1}{n}\right) = \sum_{n=1}^{\infty} \frac{c}{c-1} (0.5\sqrt{n} + 1) \cdot \left(\frac{n-1-n}{(n+1)n}\right) =$$

$$= -0.5 \cdot \frac{c}{c-1} \sum_{k=1}^{\infty} \frac{\sqrt{n}}{(n+1) \cdot n} - \frac{c}{c-1} \cdot \sum_{k=1}^{\infty} \frac{1}{n^2 + n} = -2 \cdot \frac{c}{c-1}$$
, where

$c = \{1,1;1,2;1,3\}.$

The value of the function $\varsigma(1) = \frac{1}{k_1}$ is situated in the range $\frac{c-1}{-2c} < \varsigma(1) < 1$.

Then we write down that if $s = 1$

$$L(E, s) = \frac{\varsigma(s) \cdot \varsigma(1-s)}{\varsigma(E, s)} = \frac{-0.5 \cdot \prod_p (1-p^{-s}) \cdot (1-p^{-1-s})}{\prod_p (1 - a_p p^{-s} + p^{1-2s})}, \tag{13}$$

And if $s = 1 + \varepsilon$ and when $1 < k_1 < \frac{c-1}{-2 \cdot c}$ we will have

$$L(E,s) = -0.5 \cdot \prod_3 \frac{(pp^\varepsilon - 1) \cdot (p^\varepsilon - 1)}{(pp^\varepsilon - a_p p^\varepsilon + 1)};$$

Because from lemma 1 [9, p. 33]: ' When $\operatorname{Re} s > 0$, $N >> 1$

$$\varsigma(s) = \sum_{n=1}^{N} \frac{1}{n^s} + \frac{N^{1-s}}{s-1} - \frac{1}{2} \cdot N^{-s} + s \cdot \int_N^\infty \frac{\frac{1}{2} - \{u\}}{u^{s+1}} du \text{ " we will deduce } \varsigma(0) = -0.5; \text{ when } \varepsilon \to 0$$

we will deduce

$$L(E,s) = -0.5 \cdot \lim_{\varepsilon \to 0} \prod_p \left( \frac{pp^\varepsilon - 1}{pp^\varepsilon + 1 - a_p \cdot p^\varepsilon} \right) \cdot (p^\varepsilon - 1) = -0.5 \cdot \prod_p \left( \frac{p-1}{p+1-a_p} \right) \cdot \lim_{\varepsilon \to 0} (p^\varepsilon - 1) = 0 \cdot$$

Because when $\varepsilon \to 0$ $(p^\varepsilon - 1) \to 0$, and $\left( \frac{p-1}{p+1-a_p} \right) > 1$ and if $p \to \infty$ $\lim_{p \to \infty} \frac{p-1}{p+1-a_p} = 1$

The theorem is proved.

2. Order of a zero. When the function $L(E.s)$, which is identically not equal to zero, is holomorphic in domain $D$, and equal to zero in point $a$ of this domain, then its decomposition for somewhat domain of the point is the following

$$L(E,s) = c_1 \cdot (s-1) + c_2 \cdot (s-1)^2 + ... + c_n \cdot (s-1)^n + ..., \qquad (14)$$

by virtue of the fact that $c_0 = L(E,1) = 0$.

Obviously, all coefficients $c_n$ of decomposition (14) cannot be equal to zero, by virtue of the fact that the function $L(E,s)$ is equal to zero everywhere in some neighborhood of the point $a$, by the theorem of unique solution will be identical zero in domain $D$. So, among coefficients $c_n$ $(n = 1,2,3,...)$ we have different from zero; let us denote by $n$, if $n >> 1$ -the lowest number of this coefficients.

Then we will get:

$$c_1 = c_2 = ... = c_{n-1} = 0, \qquad c_n \neq 0.$$

Now the decomposition (14) will look as follows:

$$L(E,s) = c_n \cdot (s-1)^n + c_{n+1}(s-1)^{n+1} + ..., \qquad (15)$$

where $c_n \neq 0$.

In this case point $a$ will be zero in order m for $L(E,s)$ function in point $s = 1$.

**3. Constructing the Galois group of an elliptic curve.** Let the elliptic curve $E$ be defined over the field $K$ and let the $L$ be the expansion of field $K$. It denotes following, $\sigma$ is the isomorphism of the field $L$, not certainly identical to $K$ [4, p. 27]. It defines the curve $E^\sigma$, obtained using $\sigma$ to equation coefficients, which set the curve $E$. For example, when the curve $E$ is set by equation

$$y^2 = x^3 - Ax - B,$$

Than $E^\sigma$ is defined by equation

$$y^2 = x^3 - A^\sigma x - B^\sigma.$$

When $P$, $Q$ are points on the curve $E$ in domain $L$, we can have the formula

$$(P+Q)^\sigma = P^\sigma + Q^\sigma.$$

The Total in by left side is relative to the addition to $E$, and the total in right side is related to addition to $E^\sigma$. Obviously the equality follows from the fact that the algebraic addition formula is given by rational functions of the coordinates with coefficients from the field $K$. Also, when $P = (x, y)$, than $P^\sigma = (x^\sigma, y^\sigma)$ is obtained using $\sigma$ to coordinates.

Especially, let us assume that $P$ is a point of finite order, because $NP = 0$. Since the point $O$ is rational on $K$, then any isomorphism $\sigma$ of the field $L$ over $K$ we get $NP^\sigma = 0$ and, so, $P^\sigma$ also is a point of order $N$. Then as the number of points of order $N$ is finite, it follows that they are algebraic over $K$ (that is their coordinates are algebraic over $K$).

When $P = (x, y)$, we set $K(P) = K(x, y)$ field expansion $K$, obtained by addition of coordinates of point $P$, Similarly, $k(E_N)$ the denote the composite the fields $K(P)$ for all $P \in E_n$. We emphasize that we consider all points of finite order, as points with coordinates with fixed algebraic closure of the field $K$, which we denote as $^a K$ or $K_a$.

The above remark shows that the Galois group Gal $(K_\sigma / K)$ acts as a group of elements of the set $E_N$. So, $K(E_N)$ is a normal expansion of the field and $K$ is the Galois expansion, if $N$ is not divisible by the characteristics of the field $K$. We set $K(E_N)$ as the field of points in the order $N$ of the curve $A$ over the field $K$.

Besides when $\sigma$ is automorphism of the field $K(E_N)$ over $K$ and if $\{t_1\}$, $\{t_1, t_2\}$,..., $\{t_1, t_2, ..., t_r\}$–bases $E_N$ over $Z/nZ$, than $\sigma$ can be set as matrixes

$$(a), \quad \begin{pmatrix} a_{1,1} & a_{1,2} \\ a_{2,1} & a_{2,2} \end{pmatrix}, \quad ..., \quad \begin{pmatrix} a_{1,1} & ... & a_{1,r} \\ & & \\ a_{r,1} & ... & a_{r,r} \end{pmatrix} \quad \text{such as}$$

$$(\sigma \cdot t_1), \quad \begin{pmatrix} \sigma \cdot t_1 \\ \sigma \cdot t_2 \end{pmatrix} = \begin{pmatrix} a_{1,1} t_1 + a_{1,2} t_2 \\ a_{2,1} t_1 + a_{2,2} t_2 \end{pmatrix} \cdot \begin{pmatrix} t_1 \\ t_2 \end{pmatrix} \quad \text{and etc.}$$

So we have obtained injective homomorphism:

$$Gal(K(E_N)/K) \to GL(n, Q).$$

**Theorem 2.** (**Mordell**) [1, p. 367]. Let $E$ be some elliptical curve, defined over $Q$. Then $E(Q)$ is be finitely generated abelian group.

**Theorem 3.** [1, cm. 368]. Let $E$ be some elliptical curve, defined over $Q$. Then $E(Q)$ is isomorphic to one of following groups $Z/mZ$ if $m \leq 10$ or $m = 12$, $Z/2Z \oplus Z/2mZ$ if $m \leq 4$.

**4. Theorem 4 (consistency between the group and the rank and the order of zero).** Let $L/K$ is the finite Galois expansion of degree $n$ and $\sigma_1, \sigma_2, ..., \sigma_n$ – are elements of the set $G$, where $\sigma_1 = (s-1), \sigma_2 = (s-1)^2, ..., \sigma_n = (s-1)^n$. Then there is element $\omega \in L$, such as, $\sigma_1 \cdot \omega, \sigma_2 \cdot \omega, ..., \sigma_n \cdot \omega$ create basis $L$ over $K$, then the elements of the Galois group use up ($n-1$) first coefficients of the series (14) into zeroes of series (15). So the rank of the Galois group will equal to the number of zeroes of order. $L(E,.s)$.

**Proof.** For any $\sigma \in G$ let the $X_\sigma$ be the variable and $t_{\sigma,\tau} = X_{\sigma^{-1}\tau}$. We set

$$X_i = X_{\sigma_i}. \text{ Where } \quad X_i = (s-1)_{\sigma_i}, \text{ and } \quad t_{\sigma\tau} = (s-1)_{\sigma^{-1}\tau}.$$

Let the $f(x_1, x_2, ..., x_n) = \det(t\sigma_i, \sigma_j)$.

Then $f$ is not identical to zero, that is evident, in the theorem 19 [2, p. 259] the determinant can not be equal to zero with all $x \in L$, when we in $f$ substitute $\sigma_i(x)$ instead of $X_i$. That's why there exists the element $\omega \in L$, which is

$$\det(\sigma_i^{-1} \cdot \sigma_j(\omega)) \neq 0.$$

We set coefficients of a power series (14) by $c_1, c_2, ..., c_{n-1}$. And let us assume that the elements (coefficients of power series) $c_1, c_2, ..., c_{n-1} \in K$ such as

$$c_1 \cdot \sigma_1(\omega) + c_2 \cdot \sigma_2(\omega) + ... + c_{n-1} \cdot \sigma_n(\omega) = 0.$$

We use $\sigma_i^{-1}$ according to the expression for each $i = 1, 2, ..., n-1$. Since $c_{i,j} \in K$, we get a system of simple linear equations for unknown quantity $c_j$ and receive that, $c_j = 0$ for $i = 1, 2, ..., n-1$. And so, $\omega$ will be the sought for element, in this case for $\omega = \dfrac{1}{s-1}$.

According to the Corollary of Lemma 2.3 [10 st.144]: "Let the $L-$ be the finite expansion of the field $K$ with Abelian Galois group $G$ of power which divides $n$. Then the group $G$ is a direct product of cyclic subgroups $G_1, G_2, ..., G_r$. Suppose that for each $i$ over $L_i$ will be set the subfield, fixed for the subgroups $G_1 \times G_2 \times ... \times G_r$; then $G(L_i / K) = G_i$, $L_i = K(\alpha_i)$, where $\alpha_i^n = a_i \in K$ and $L = K(\alpha_1,...,\alpha_n)$".

And Lemmas 2.4[10, p.144]. When $L$ is the normal algebraic expansion of $K$ with the Galois group $G$, then

$$H^1(G, L^*) = 0$$

Then using Theorem 1 and Lemma 2.4, we can write that in the normal expansion of the Galois group, change series (14) into series

$$L(E,1) = c_n (s-1)^n + c_{n+1}(s-1)^{n+1} + .... \qquad (16)$$

Using the formula (16) we will have, if $n = r_E$

$$C_{n,E} = \lim_{s \to 1+0} \frac{L(E,1)}{(s-1)^{r_E}} = \lim_{s \to 1+0} \frac{c_n(s-1)^n + c_{n+1}(s-1)^{n+1}}{(s-1)^{r_E}} = c_n.$$

According to Lemma 1 and the result [9, p.33] the function $L(E,s)$ is analytical in whole domain $s \in (0, \infty)$, it can be expanded into a Taylor series in powers of $(s-1)$ and with coefficients $C_{n,E} = \frac{L^{(n)}(E,1)}{k!} = \frac{B_{n,E}}{k!}$, where $L^{(n)}(E,s)$ –is derivative by the order $n$ – from power of Hasse-Weil function $L(E,s)$. Thus, the rank order of the Galois group is equal to zero Hasse-Weil function.

**Theorem is proved.** So the Birch and Swinnerton-Dyer conjecture is fair.

**Key words**: the hypothesis of Birch and Swinnerton-Dyer, function of Hasse-Weil , Riemann's hypothesis, the Galois group, the complex power series.

05 September 2013

Rewriten 14 May 2014